\documentclass[11pt,leqno,a4paper]{article} 
\usepackage{graphics}
\newtheorem{thm}{Theorem}[section]

\newtheorem{prop}{Proposition}

\newcommand{\beqa}{\begin{eqnarray}}
\newcommand{\eeqa}{\end{eqnarray}}

\newcommand{\pf}{\noindent {\bf Proof:} $\s$ }
\newcommand{\epf}{ \hfill$\diamondsuit$ \medskip}
\newcommand{\md}{\medskip}

\newcommand{\beq}{\begin{equation}}
\newcommand{\eeq}{\end{equation}}
\newcommand{\lbl}{\label}
\newcommand{\s}{\; \;}

\newcommand{\al}{\alpha}

\title{Non-existence of solutions for non-autonomous elliptic systems}

\author{
Philip Korman  \thanks{Supported in part by the Taft Faculty Grant at
the University of Cincinnati} \\ 
Department of Mathematical Sciences \\ 
University of Cincinnati \\ 
Cincinnati Ohio 45221-0025 \\
}

\date{}

\begin{document}

\maketitle
\begin{abstract} 
We extend the classical Pohozaev's identity to semilinear elliptic systems  of Hamiltonian type, providing a simpler approach, and a generalization,  of the results  of E. Mitidieri \cite{M}, R.C.A.M. Van der Vorst \cite{V}, and Y. Bozhkov and E. Mitidieri \cite{BM}.
 \end{abstract}

\begin{flushleft}
\end{flushleft}

\begin{flushleft}
\end{flushleft}

\section{Introduction}
\setcounter{equation}{0}
\setcounter{thm}{0}
\setcounter{lma}{0}

Any solution $u(x)$ of semilinear Dirichlet problem on a bounded domain $\Omega \subset R^n$
\beq
\lbl{1}
\Delta u +f(x,u)=0 \s \mbox{in $\Omega$}, \s u=0 \s \mbox{on $ \partial \Omega$}
\eeq
satisfies the well known Pohozaev's identity
\beq
\lbl{2}
\int _{\Omega} \left[2n F(x,u)+(2-n)uf(x,u) +2\Sigma _{i=1}^n x_iF_{x_i}(x,u) \right] \, dx=\int _{\partial \Omega} (x \cdot \nu) | \nabla u |^2 \, dS \,.
\eeq
Here $F(x,u)=\int_0^u f(x,t) \,dt$, and $\nu$ is the unit normal vector on $\partial \Omega$, pointing outside. (From the equation (\ref{1}), $
\int _{\Omega} uf(x,u) \, dx= \int _{\Omega} |\nabla u |^2 \, dx$, which gives an alternative form of the Pohozaev's identity.) Pohozaev's identity is usually written for the  case $f=f(u)$, but the present version is also known, see e.g., K. Schmitt \cite{S}. A standard use of this identity is to conclude that if $\Omega$ is a star-shaped domain with respect to the origin, i.e., $x \cdot \nu \geq 0$ for all 
$x \in  \partial \Omega$, and $f(u)=u|u|^{p-1}$, for some constant $p$, then the problem (\ref{1}) has no non-trivial solutions in the super-critical case, when $p>\frac{n+2}{n-2}$. In this note we present a proof of Pohozaev's identity, which appears a little more straightforward than the usual one, see e.g., L. Evans \cite{E}, and then use a similar idea for systems, generalizing  the well-known results of E. Mitidieri \cite{M}, see also R.C.A.M. Van der Vorst \cite{V}, and of   Y. Bozhkov and E. Mitidieri \cite{BM}, by allowing explicit dependence on $x$ in the Hamiltonian function. 
\md 

Let $z=x \cdot  \nabla u= \Sigma _{i=1}^n x_i u_{x_i}$. It is straightforward  to verify that $z$ satisfies 
\beq
\lbl{3}
\Delta z +f_u(x,u)z=-2f(x,u)-\Sigma _{i=1}^n x_i f_{x_i}(x,u) \,.
\eeq
We multiply the equation (\ref{1}) by $z$, and subtract from that the equation (\ref{3}) multiplied by $u$, obtaining 
\beq
\lbl{4}
\Sigma _{i=1}^n  \left(zu_{x_i}-uz _{x_i} \right)_{x_i} +\Sigma _{i=1}^n \left(f(x,u)-uf_u(x,u) \right) x_i u_{x_i}=
2f(x,u)u+\Sigma _{i=1}^n x_i f_{x_i}(x,u) u \,.
\eeq
We have 
\beqa \nonumber
& \Sigma _{i=1}^n \left(f(x,u)-uf_u(x,u) \right) x_i u_{x_i}=\Sigma _{i=1}^n  x_i \frac{\partial}{\partial x_i} (2F-uf)-2
\Sigma _{i=1}^n x_iF_{x_i}+\Sigma _{i=1}^n x_i f_{x_i}(x,u) u= \\ \nonumber
& \Sigma _{i=1}^n   \frac{\partial}{\partial x_i} \left[ x_i(2F-uf) \right]-n(2F-uf) -2
\Sigma _{i=1}^n x_iF_{x_i}+\Sigma _{i=1}^n x_i f_{x_i}(x,u) u \,. \nonumber
\eeqa
We then rewrite (\ref{4})
\beq
\lbl{5}
\s\s\s\s \Sigma _{i=1}^n \left[ (zu_{x_i}-uz _{x_i}) +x_i  (2F(x,u)-uf(x,u))\right]_{x_i} =2n F(x,u)+(2-n)uf(x,u)+2\Sigma _{i=1}^n x_iF_{x_i} \,.
\eeq
Integrating over $\Omega$, we conclude the Pohozaev's identity (\ref{2}). (The only non-zero boundary term is $\Sigma _{i=1}^n  \int _{\partial \Omega} z u_{x_i} \nu _i \, dS$. Since $\partial \Omega$ is a level set of $u$, $\nu =\pm \frac{\nabla u}{|\nabla u|}\;$, i.e., 
$u_{x_i}=\pm |\nabla u| \nu _i$. Then $z=\pm (x \cdot \nu) | \nabla u |$, and $\Sigma _{i=1}^n u_{x_i} \nu _i=\pm | \nabla u |$.)
\md 

We refer to (\ref{5}) as a {\em differential form} of Pohozaev's identity. For radial solutions on a ball, the corresponding version of (\ref{5}) played a crucial role in the study of exact multiplicity of solutions, see T. Ouyang and J. Shi \cite{OS1}, and also P. Korman \cite{K40}, which shows the potential usefulness of this identity.

\section{Non-existence of solutions for a class of systems}
\setcounter{equation}{0}
\setcounter{thm}{0}
\setcounter{lma}{0}

The following class of systems has attracted considerable attention recently
\beqa
\lbl{6}
& \Delta u +H_v(u,v)=0 \s \mbox{in $\Omega$}, \s u=0 \s \mbox{on $ \partial \Omega$} \\
& \Delta v +H_u(u,v)=0 \s \mbox{in $\Omega$}, \s v=0 \s \mbox{on $ \partial \Omega$} \,, \nonumber
\eeqa
where $H(u,v)$ is a given differentiable function, 
see e.g., the following surveys: D.G. de Figueiredo \cite{F}, P. Quittner and P.  Souplet \cite{Q}, B. Ruf \cite{R}, see also P. Korman \cite{K26}. This system is of {\em Hamiltonian} type, so that it has some of the properties of scalar equations.  
\md 

More generally, let $H=H(x,u_1,u_2, \ldots, u_m,v_1,v_2, \ldots, v_m)$, with  integer $m \geq 1$, and consider the Hamiltonian system of $2m$ equations
\beqa
\lbl{6.1}
& \Delta u_k +H_{v_k}=0 \s \mbox{in $\Omega$}, \s u_k=0 \s \mbox{on $ \partial \Omega$}, \s k=1,2, \ldots, m \\
& \Delta v_k +H_{u_k}=0 \s \mbox{in $\Omega$}, \s v_k=0 \s \mbox{on $ \partial \Omega$}, \s k=1,2, \ldots, m \,. \nonumber
\eeqa

We call solution of (\ref{6.1}) to be positive, if $u_k(x)>0$ and $v_k(x)>0$ for all $x \in \Omega$, and all $k$.
We consider only the classical solutions, with $u_k$ and $v_k$  of class $C^2(\Omega) \cap C^1(\bar \Omega)$. We have the following generalization of the results of  \cite{BM} and \cite{M}.

\begin{thm}\lbl{thm:1}
Assume that $H(x,u_1,u_2, \ldots, u_m,v_1,v_2, \ldots, v_m) \in C^2( \Omega \times R^m_+  \times  R^m_+ ) \cap C(\bar \Omega \times \bar R^m_ + \times \bar R^m_+ )$ satisfies
\beq
\lbl{6.2}
H(x,0, \ldots, 0,0, \ldots, 0)=0 \s \mbox{for all $x \in \partial \Omega$} \,.
\eeq
Then
 for any positive solution of (\ref{6.1}), and any real numbers  $a_1, \ldots, a_m$, one has
\beqa
\lbl{7}
&  \int _{\Omega} \left[2n H +(2-n) \Sigma _{k=1}^m \left(  a_k u_kH_{u_k}+ (2-a_k) v_kH_{v_k} \right) +2 \Sigma _{i=1}^n  x_i H_{x_i} \right] \, dx  \\
& =2 \Sigma _{k=1}^m \int _{\partial \Omega} (x \cdot \nu) | \nabla u_k | | \nabla v_k |\, dS \,. \nonumber
\eeqa
\end{thm}

\pf
Define $p_k=x \cdot  \nabla u_k= \Sigma _{i=1}^n x_i u_{kx_i}$, and $q_k=x \cdot  \nabla v= \Sigma _{i=1}^n x_i v_{kx_i}$, $ k=1,2, \ldots, m$. These functions satisfy the system
\beqa
\lbl{8}
&  \Delta p_k +\Sigma _{j=1}^m H_{v_ku_j}p_j+\Sigma _{j=1}^m H_{v_kv_j}q_j=-2H_{v_k}-\Sigma _{i=1}^n x_i H_{v_kx_i} ,  \s k=1,2, \ldots, m\\
& \Delta q_k +\Sigma _{j=1}^m H_{u_ku_j}p_j+\Sigma _{j=1}^m H_{u_k v_j}q_j=-2H_{u_k}-\Sigma _{i=1}^n x_i H_{u_kx_i},  \s k=1,2, \ldots, m \,. \nonumber
\eeqa
We multiply the first equation in (\ref{6.1}) by $q_k$, and subtract from that the first equation in (\ref{8}) multiplied by $v_k$. The result can be written as
\beqa
\lbl{9}
&  \Sigma _{i=1}^n \left[ (u_{kx_i}q_k-p _{kx_i}v_k)_{x_i} +( -u_{kx_i}q_{kx_i}+v _{kx_i}p_{kx_i}) \right]  \\
& +H_{v_k}q_k-\Sigma _{j=1}^m H_{v_ku_j}p_jv_k-\Sigma _{j=1}^m H_{v_kv_j}q_jv_k=2v_kH_{v_k} +v_k \Sigma _{i=1}^n x_i H_{v_kx_i} \,. \nonumber
\eeqa
Similarly, we multiply the second  equation in (\ref{6.1}) by $p_k$, and subtract from that the second  equation in (\ref{8}) multiplied by $u_k$, and write the result  as
\beqa
\lbl{10}
&  \Sigma _{i=1}^n \left[ (v_{kx_i}p_k-q _{kx_i}u_k )_{x_i} +( -v_{kx_i}p_{kx_i}+u _{kx_i}q_{kx_i}) \right]  \\
& +H_{u_k}p_k-\Sigma _{j=1}^m H_{u_ku_j}p_ju_k-\Sigma _{j=1}^m H_{u_kv_j}q_ju_k=2u_kH_{u_k}+u_k \Sigma _{i=1}^n x_i H_{u_kx_i} \,. \nonumber
\eeqa
Adding the equations (\ref{9}) and (\ref{10}), we get 
\beqa  \nonumber
&   \Sigma _{i=1}^n \left[ u_{kx_i}q_k-p _{kx_i}v_k+ v_{kx_i}p_k-q _{kx_i}u_k \right]_{x_i} +H_{u_k}p_k +H_{v_k}q_k-\Sigma _{j=1}^m H_{u_ku_j}p_ju_k\\ \nonumber
& -\Sigma _{j=1}^m H_{u_kv_j}q_ju_k-\Sigma _{j=1}^m H_{v_ku_j}p_jv_k-\Sigma _{j=1}^m H_{v_kv_j}q_jv_k \\
& =2u_kH_{u_k}+2v_kH_{v_k}+ u_k \Sigma _{i=1}^n x_i H_{u_kx_i} +v_k \Sigma _{i=1}^n x_i H_{v_kx_i} \,. \nonumber
\eeqa
We now sum in $k$,   putting  the result into the form
\beqa  \nonumber
&  \Sigma _{k=1}^m \Sigma _{i=1}^n \left[ u_{kx_i}q_k-p _{kx_i}v_k+ v_{kx_i}p_k-q _{kx_i}u_k \right]_{x_i}  \\ \nonumber
&  +\Sigma _{i=1}^n x_i  \left(2H-\Sigma _{k=1}^m u_kH_{u_k}-\Sigma _{k=1}^m v_kH_{v_k} \right) _{x_i}=2\Sigma _{k=1}^m u_kH_{u_k}+2 \Sigma _{k=1}^m v_kH_{v_k} +2 \Sigma _{i=1}^n  x_i H_{x_i}\,. \nonumber
\eeqa
Writing, 
\beqa  \nonumber
& \Sigma _{i=1}^n  x_i \frac{\partial}{\partial x_i} (2H-\Sigma _{k=1}^m u_kH_{u_k}-\Sigma _{k=1}^m v_kH_{v_k})=\Sigma _{i=1}^n   \frac{\partial}{\partial x_i} \left[ x_i(2H-\Sigma _{k=1}^m u_kH_{u_k}-\Sigma _{k=1}^m v_kH_{v_k}) \right] \\ \nonumber
& -n(2H-\Sigma _{k=1}^m u_kH_{u_k}-\Sigma _{k=1}^m v_kH_{v_k}) \,, \nonumber
\eeqa
we obtain  the differential form of Pohozaev's identity
\beqa  \nonumber
&   \Sigma _{k=1}^m \Sigma _{i=1}^n \left[ u_{kx_i}q_k-p _{kx_i}v_k+ v_{kx_i}p_k-q _{kx_i}u_k  +x_i  \left(2H-\Sigma _{k=1}^m u_kH_{u_k}-\Sigma _{k=1}^m v_kH_{v_k} \right) \right]_{x_i} \\ \nonumber
&  =2n H+(2-n) \left(\Sigma _{k=1}^m u_kH_{u_k}+ \Sigma _{k=1}^m v_kH_{v_k} \right)+2 \Sigma _{i=1}^n  x_i H_{x_i} \,. \nonumber
\eeqa
Integrating, we obtain, in view of (\ref{6.2}),
\beqa
\lbl{11}
&  \int _{\Omega} \left[2n H(u,v) +(2-n)\left( \Sigma _{k=1}^m u_kH_{u_k}+ \Sigma _{k=1}^m v_kH_{v_k} \right) +2 \Sigma _{i=1}^n  x_i H_{x_i} \right] \, dx  \\
& =2 \Sigma _{k=1}^m \int _{\partial \Omega} (x \cdot \nu) | \nabla u_k | | \nabla v_k |\, dS \,. \nonumber
\eeqa
(Since we consider positive solutions, and $\partial \Omega$ is a level set for both $u_k$ and $v_k$, we have $\nu =- \frac{\nabla u_k}{|\nabla u_k|}=- \frac{\nabla v_k}{|\nabla v_k|}$, i.e., $u_{ki}=-|\nabla u_k| \nu _i$ and $v_{ki}=-|\nabla v_k| \nu _i$ on the boundary $\partial \Omega$.)
From the first equation in (\ref{6.1}), $\int _{\Omega} v_kH_{v_k} \, dx=\int _{\Omega} \nabla u_k \cdot \nabla v_k \, dx$, while from the second equation $\int _{\Omega} u_k H_{u_k} \, dx=\int _{\Omega} \nabla u_k \cdot \nabla v_k \, dx$, i.e., for each $k$
\[
\int _{\Omega} v_k H_{v_k} \, dx=\int _{\Omega} u_k H_{u_k} \, dx \,.
\]
Using this in (\ref{11}), we conclude the proof.
\epf

\noindent
{\bf Remarks} 
\begin{enumerate}
  \item 
We consider only the classical solutions. Observe that by our conditions and elliptic regularity, classical solutions are in fact of class $C^3(\Omega)$, so that all quantities in the above proof are well defined. 
\item
In case $H$ is independent of $x$, the condition (\ref{6.2}) can be assumed without loss of generality.
\end{enumerate}
\md 

As a consequence, we have the following non-existence result.

\begin{prop}
Assume that $\Omega$ is a star-shaped domain with respect to the origin, and  for some real constants $\al _1, \ldots, \al _m$,  all $u_k>0$, $v_k>0$, and all $x \in \Omega$,  we have
\beq
\lbl{100.1}
nH+(2-n)  \Sigma _{k=1}^m \left( \al _k u_kH_{u_k}+(1-\al _k) v_kH_{v_k} \right)+ \Sigma _{i=1}^n  x_i H_{x_i}  <0 \,.
\eeq
Then the problem (\ref{6.1}) has no positive solutions.
\end{prop}

\pf
We use the identity (\ref{7}), with $a_k/2=\al _k$. Then, assuming existence of positive solution, the left hand side of (\ref{7}) is negative, while the right hand side is non-negative, a contradiction.
\epf

Observe, that it suffices to assume that $\Omega$ is star-shaped with respect to any one of its points (which we then take to be the origin).
\md

In case $m=1$, and $H=H(u,v)$, we recover the following condition of E. Mitidieri \cite{M}. 

\begin{prop}
Assume that $\Omega$ is a star-shaped domain with respect to the origin, and  for some real constant $\al$, and all $u>0$, $v>0$  we have
\beq
\lbl{10.1}
 \al uH_u(u,v)+ (1-\al) vH_v(u,v)> \frac{n}{n-2} H(u,v) \,.
\eeq
Then the problem (\ref{6}) has no positive solutions.
\end{prop}

Comparing this result to E. Mitidieri \cite{M}, observe that we do not require that $H_u(0,0)=H_v(0,0)=0$.
\md

An important subclass of (\ref{6}) is
\beqa
\lbl{11.1}
& \Delta u +f(v)=0 \s \mbox{in $\Omega$}, \s u=0 \s \mbox{on $ \partial \Omega$} \\
& \Delta v +g(u)=0 \s \mbox{in $\Omega$}, \s v=0 \s \mbox{on $ \partial \Omega$} \,, \nonumber
\eeqa
which corresponds to $H(u,v)=F(v)+G(u)$, where  $F(v)=\int_0^v f(t) \,dt$, $G(u)=\int_0^u g(t) \,dt$. Unlike \cite{M}, we do not require that $f(0)=g(0)=0$. 
The Theorem \ref{thm:1} now reads as follows.

\begin{thm}\lbl{thm:3}
Let $f, \, g \in C(\bar R_+)$.
For any positive solution of (\ref{11.1}), and any real number  $a$, one has
\beqa
\lbl{11.2}
&  \int _{\Omega} \left[2n (F(v)+G(u))+(2-n)\left(a vf(v)+(2-a)ug(u)  \right) \right] \, dx  \\
& =2 \int _{\partial \Omega} (x \cdot \nu) | \nabla u | | \nabla v |\, dS \,. \nonumber
\eeqa
\end{thm}
\md 

More generally, we consider 
\beqa
\lbl{11.1a}
& \Delta u +f(x,v)=0 \s \mbox{in $\Omega$}, \s u=0 \s \mbox{on $ \partial \Omega$} \\
& \Delta v +g(x,u)=0 \s \mbox{in $\Omega$}, \s v=0 \s \mbox{on $ \partial \Omega$} \,, \nonumber
\eeqa
with $H(x,u,v)=F(x,v)+G(x,u)$, where  $F(x,v)=\int_0^v f(x,t) \,dt$, $G(x,u)=\int_0^u g(x,t) \,dt$. 

\begin{thm}\lbl{thm:3a}
Let $f, \, g \in C(\Omega \times \bar R_+)$.
For any positive solution of (\ref{11.1a}), and any real number  $a$, one has
\beqa
\lbl{11.2a}
&  \s\s\s \int _{\Omega} \left[2n (F(x,v)+G(x,u))+(2-n)\left(a vf(x,v)+(2-a)ug(x,u) \right) +2\Sigma _{i=1}^n  x_i \left(  F_{x_i} +G_{x_i } \right)  \right] \, dx  \\
& =2 \int _{\partial \Omega} (x \cdot \nu) | \nabla u | | \nabla v |\, dS \,. \nonumber
\eeqa
\end{thm}
We now consider a particular system
\beqa
\lbl{12}
& \Delta u +v^p=0 \s \mbox{in $\Omega$}, \s u=0 \s \mbox{on $ \partial \Omega$} \\
& \Delta v +g(x,u)=0 \s \mbox{in $\Omega$}, \s v=0 \s \mbox{on $ \partial \Omega$} \,, \nonumber
\eeqa
with $g(x,u) \in C(\Omega \times \bar R_+)$, and a constant $p>0$.

\begin{thm}\lbl{thm:2}
Assume that $\Omega$ is a star-shaped domain with respect to the origin,   and 
\beq
\lbl{15}
\s\s\s\s n G(x,u)+(2-n) \left(1-\frac{n}{(n-2)(p+1)} \right)ug(x,u)+\Sigma _{i=1}^n  x_i  G_{x_i }<0 \,,  \s \mbox{for $x \in \Omega$, and  $u>0$} \,.
\eeq
Then the problem (\ref{12}) has no positive solutions.
\end{thm}

\pf
We use  the identity (\ref{11.2a}), with $f(v)=v^p$. We select the constant $a$, so that
\[
2n F(v)+(2-n)a vf(v)=0 \,,
\]
i.e., $a=\frac{2n}{(n-2)(p+1)}$. Then, assuming existence of a positive solution, the left hand side of (\ref{11.2a}) is negative, while the right hand side is non-negative, a contradiction.
\epf

Observe that in case $p=1$, the Theorem \ref{thm:2} provides a non-existence result for a biharmonic problem with Navier boundary conditions
\beq
\lbl{14}
\Delta ^2 u=g(x,u) \s \mbox{in $\Omega$}, \s u=\Delta u=0 \s \mbox{on $ \partial \Omega$} \,.
\eeq

\begin{prop} 
Assume that $\Omega$ is a star-shaped domain with respect to the origin,   and the condition (\ref{15}), with $p=1$,  holds. Then the problem (\ref{14}) has no positive solutions.
\end{prop}

Finally, we consider the  system
\beqa
\lbl{16}
& \Delta u +v^p=0 \s \mbox{in $\Omega$}, \s u=0 \s \mbox{on $ \partial \Omega$} \\
& \Delta v +u^q=0 \s \mbox{in $\Omega$}, \s v=0 \s \mbox{on $ \partial \Omega$} \,. \nonumber
\eeqa

The curve $\frac{1}{p+1}+\frac{1}{q+1}=\frac{n-2}{n}$ is called a {\em critical hyperbola}. We recover the following well known result of E. Mitidieri \cite{M}, see also R.C.A.M. Van der Vorst \cite{V}. (Observe that we relax the restriction $p$, $q>1$ from \cite{M}.)

\begin{prop} 
Assume that $p$, $q>0$, and 
\beq
\lbl{17}
\frac{1}{p+1}+\frac{1}{q+1}< \frac{n-2}{n} \,.
\eeq
Then the problem (\ref{16}) has no positive solutions.
\end{prop}

\pf
Condition (\ref{17}) implies  (\ref{15}), and then the Theorem \ref{thm:2} applies.
\epf

In case $p=1$, we recover the following known result, see E. Mitidieri \cite{M}.
\begin{prop}
Assume that $\Omega$ is a star-shaped domain with respect to the origin, and $q>\frac{n+4}{n-4}$.  Then the problem 
\beq
\lbl{18}
\Delta ^2 u=u^q \s \mbox{in $\Omega$}, \s u=\Delta u=0 \s \mbox{on $ \partial \Omega$} 
\eeq
has no positive solutions.
\end{prop}


\begin{thebibliography}{99}


\bibitem{BM}
Y. Bozhkov and E. Mitidieri, The Noether approach to Pokhozhaev's identities, {\em  Mediterr. J. Math.}  {\bf 4},  no. 4, 383-405   (2007).
\vspace{-0.2cm}

\bibitem{E}
L. Evans, Partial Differential Equations. Graduate Studies in Mathematics, 19. American Mathematical Society, Providence, RI, 1998.
\vspace{-0.2cm}

\bibitem{F}
D.G. de Figueiredo, Semilinear elliptic systems: existence, multiplicity, symmetry of solutions, Handbook
of Differential Equations, Stationary Partial  Differential Equations,  Vol. {\bf 5}, Edited by M. Chipot, Elsevier Science, North Holland, 1-48 (2008).
\vspace{-0.2cm}

\bibitem{K26}
P. Korman,    Pohozaev's identity and non-existence of solutions for elliptic systems, {\em  Comm. Appl. Nonlinear Anal.} {\bf  17}, no. 4, 81-88  (2010).
\vspace{-0.2cm}

\bibitem{K40}
P. Korman, Uniqueness and exact multiplicity of solutions for non-autonomous Dirichlet problems, {\em  Adv. Nonlinear Stud.} {\bf   6},  no. 3, 461-481   (2006).
\vspace{-0.2cm}


\bibitem{M}
E. Mitidieri,  A Rellich type identity and applications, {\em   Comm. Partial Differential Equations} {\bf  18},  no. 1-2, 125-151   (1993). 
\vspace{-0.2cm}

\bibitem{OS1}
T. Ouyang and J. Shi,
Exact multiplicity of positive solutions for a class of semilinear problems, II, {\em J. Differential Equations} {\bf 158},  no. 1, 94-151 (1999).
\vspace{-0.2cm}

\bibitem{P}
S. I. Pohozaev,  On the eigenfunctions of the equation $\Delta u+\lambda f(u)=0$. (Russian)  {\em Dokl. Akad. Nauk SSSR}  {\bf 165}, 36-39   (1965). 
\vspace{-0.2cm}

\bibitem{PuS1}
P. Pucci and J.  Serrin,  A general variational identity, {\em Indiana Univ. Math. J.} {\bf 35},  no. 3, 681-703   (1986).
\vspace{-0.2cm}

\bibitem{R}
F. Rellich,  Darstellung der Eigenwerte von $\Delta u+\lambda u=0$ durch ein Randintegral. (German)  {\em Math. Z.} {\bf   46}, 635-636  (1940).
\vspace{-0.2cm}

\bibitem{Q}
P. Quittner and P.  Souplet,  Superlinear Parabolic Problems. Blow-up, global existence and steady states. Birkh\"{a}user Advanced Texts: Basler Lehrb\"{u}cher.  Birkh\"{a}user Verlag, Basel,  (2007). 
\vspace{-0.2cm}

\bibitem{R} 
B. Ruf, Superlinear elliptic equations and systems, Handbook
of Differential Equations, Stationary Partial  Differential Equations,  Vol. {\bf 5}, Edited by M. Chipot, Elsevier Science, North Holland, 277-370. (2008).
\vspace{-0.2cm}

\bibitem{S}
K. Schmitt,  Positive solutions of semilinear elliptic boundary value problems.  Topological methods in differential equations and inclusions (Montreal, PQ, 1994),  447-500, NATO Adv. Sci. Inst. Ser. C Math. Phys. Sci., 472, Kluwer Acad. Publ., Dordrecht, (1995).
\vspace{-0.2cm}

\bibitem{V}
R.C.A.M. Van der Vorst,  Variational identities and applications to differential systems, {\em  Arch. Rational Mech. Anal.} {\bf   116},  no. 4, 375-398  (1992).
\end{thebibliography}
\end{document}